\documentstyle[11pt]{article}
\newtheorem{theo}{Theorem}[section]
\newtheorem{deff}{Definition}[section]
\newtheorem{prop}{Proposition}[section]
\newtheorem{lem}{Lemma}[section]

\newtheorem{coro}{Corollary}[section]

\input{tcilatex}

\begin{document}

\title{Fixed Point Theorems in Modular Spaces }
\date{}
\author{By \\
%EndAName
E.Hanebaly}
\maketitle

\QTP{Body Math}
\bigskip\ \ \ \ \ \ \ \ \ \ \ \ \ \ \ \ \ \ \ \ \ \ \ \ \ \ \ \ {\LARGE %
FIXED\ POINT\ THEOREMS\ in MODULAR\ SPACES}

\bigskip\ \ \ \ \ \ \ \ \ \ \ \ \ \ \ \ \ \ \ \ \ \ \ \ \ \ \ \ \ \ \ \ \ \
\ \ \ \ \ \ \ \ \ \ \ \ \ \ \ \ \ \ \ \ \ \ \ \ \ \ \ By E.Hanebaly

Abstract. By iterative techniques, we present two fixed point theorems,
whose modular formulations are relatively close to the Banach's fixed point
theorem in the normed spaces.\newline
The first result concerns the fixed point of the strongly $\rho $%
-contraction mappings. The second result deals with the fixed point of the
strict $\rho $-contraction mappings where the modular satisfies the $\Delta
_{2}$-condition.\newline
For the $\rho $-nonexpansive mappings, where the modular $\rho $ satisfies
the regular growth condition, we present a fixed point theorem of the
Schauder's type, without boundedness conditions on the domain of these
mappings.\newline
A.M.S. Subject classifications: 46A80.47H10.\newline

\section{ Introduction}

By iterative techniques, we present two fixed point theorems, whose modular
formulations are relatively close to the Banach's fixed point theorem in the
normed spaces.\newline
The first result concerns the fixed point of the strongly $\rho$-contraction
mappings.\newline
As a consequence, we get an improved version of the theorem I-1\cite{ai}, in
particular, by the deletion of the hypothesis, the $\Delta_2 $-condition and
the Fatou property.\newline
The second result concerns the fixed point of the strict $\rho$-contraction
mappings where the modular $\rho$ satisfies the $\Delta_2 $-condition. With
the last condition, the iterative techniques are to happen locally.\newline
For the $\rho$-nonexpansive mapping where $\rho$ satisfies the regular
growth condition, noted $T$ , we present one result of the Schauder's type
(i.e. ${\overline{T(B)}}^\rho $ is $\rho$-compact, with $B$ is the domain of 
$T$ ) without boundedness conditions on $B$ used in theorems 2.13 \cite{krk}
and 2.5 \cite{m1}.

\section{I- Strongly $\protect\rho $-contraction mappings}

\begin{deff}
Let $X_{\rho }$ be a modular space.\newline
. A sequence $(x_{n})_{n\in I\!\!N}$ in $X_{\rho }$ is $\rho $-convergeant
to $x\in X_{\rho }$ if: $\exists c>0$ such that $\rho
(c(x_{n}-x))\rightarrow 0$ as $n\rightarrow +\infty $.\newline
$X_{\rho }$ is $\rho $-complete if every $\rho $-Cauchy sequence $%
(x_{n})_{n\in I\!\!N}$ in $X_{\rho }$ is $\rho $-convergent, i.e. If $%
\exists c>0$ such that $\rho (c(x_{n}-x_{m}))\rightarrow 0$ as $%
n,m\rightarrow +\infty $ , then, $\exists x\in X_{\rho }$ such that $\rho
(c(x_{n}-x))\rightarrow 0$ as $n\rightarrow +\infty $.
\end{deff}

For example, Musielak-Orlicz space is $\rho $-complete in the sens of the
above definition. (cite: jm ).\newline
The following result can be considered as the first approach of the Banach's
fixed point theorem in the normed spaces.

\begin{theo}
1.1.Let $X_{\rho }$ be a $\rho $-complete modular space, and $B\subseteq
X_{\rho }$ a $\rho $-closed subset of $X_{\rho }$.\newline
Let $T:B\rightarrow B$ be a mapping such that: \newline
$\exists c,k,l\in {R\hspace*{-0.9ex}\rule{0.15ex}{1.5ex}\hspace*{0.9ex}}^{+}$
with $c>l,\ \ k\in ]0,1[$ and $\rho (c(Tx-Ty))\leq k\rho (l(x-y)),\ \
\forall x,y\in B\ \ \ (\ast )$\newline
Then $T$ has a fixed point .
\end{theo}

Remarks.\newline
We note that if $\rho (l(x-y))<+\infty \ ,\ \ \forall x,y\in B$ , then the
fixed point is unique. The insertion of the constants $c,k$ and $l$ in $%
(\ast )$ has been the field of application of this result and may be useful
(see the study, by a fixed point theorem, of an integral equation of $\rho $%
-Lipschitz or perturbed integral equations in modular function space $%
C^{\varphi }=C([0,A],L^{\varphi })\ \ \cite{ai},\cite{hh})$\newline
we note that the contraction $(\ast )$ is also valid for all constants $%
c_{0},l_{0}$ and $k_{0}$ with $l\leq l_{0}<c_{0}\leq c$ and $0<k\leq k_{0}<1$%
. Indeed: 
\begin{eqnarray*}
\rho (c_{0}(Tx-Ty)) &\leq &\rho (c(Tx-Ty)) \\
&\leq &k\rho (l(x-y)) \\
&\leq &k_{0}\rho (l_{0}(x-y))
\end{eqnarray*}
Because $\alpha \rightarrow \rho (\alpha x),\ \ (\alpha \in {\
R\hspace*{-0.9ex}\rule{0.15ex}{1.5ex}\hspace*{0.9ex}}^{+})$ is increasing.%
\newline
If $1\in \lbrack l,c]$, then $T$ is a strict $\rho $-contraction because:%
\newline
\begin{eqnarray*}
\rho (Tx-Ty) &\leq &\rho (c(Tx-Ty)) \\
&\leq &k\rho (l(x-y)) \\
&\leq &k\rho (x-y)
\end{eqnarray*}
But, as $c>l$, we have:\newline
\[
\rho (\lambda (Tx-Ty))\leq k\rho (\lambda (x-y)) 
\]
where $\lambda =l$ or $\lambda =c$ .\newline
With this last inequality, it can be said that $T$ is a strict $\rho $%
-contraction. Hence, it can be said that $T$ is a strongly $\rho $%
-contraction if $T$ satisfies $(\ast )$ in theorem 1.1. \newline
The supplementary condition $c>l$, in $(\ast )$, has permitted to delete the
boundedness condition concerning the domain of $T$ in \cite{krk}-\cite{m1}- 
\cite{Lm} where $T$ is a strict $\rho $-contraction.\newline
But, the hypothesis $c>l$ , in theorem 1.1, is essential, and that is to
apply constantly the inequality of the modular $\rho $ in the following
proof.\newline
Proof of the theorem 1.1. \newline
Let $\alpha \in {\ IR}^{+}{\hspace*{-0.9ex}}$ be the conjugate of $\frac{c}{l%
}$, i.e. $\frac{l}{c}+\frac{1}{\alpha }=1$. We assume without any loss of
generality that: \newline
$\exists x\in B$ such that $r=\rho (\alpha l(Tx-x))<+\infty $. Then the
sequence $\{T^{n}x\}_{n\in I\!\!N}$ is $\rho $-Cauchy. Indeed:\newline
\begin{eqnarray*}
\rho (c(T^{n+m}x-T^{m}x)) &\leq &k\rho (l(T^{n+m-1}x-T^{m-1}x)) \\
&\leq &k\rho (c(T^{n+m-1}x-T^{m-1}x)) \\
&\leq &k^{2}\rho (l(T^{n+m-2}x-T^{m-2}x))
\end{eqnarray*}
By induction, we deduce: 
\[
\rho (c(T^{n+m}x-T^{m}x))\leq k^{m}\rho (l(T^{n}x-x)) 
\]
Moreover, 
\begin{eqnarray*}
\rho (l(T^{n}x-x)) &=&\rho (l(T^{n}x-Tx)+l(Tx-x)) \\
&=&\rho (\frac{l}{c}c(T^{n}x-Tx)+\frac{\alpha l}{\alpha }(Tx-x)) \\
&\leq &\rho (c(T^{n}x-Tx))+\rho (\alpha l(Tx-x)) \\
&\leq &k\rho (l(T^{n-1}x-x))+r
\end{eqnarray*}
By induction, we have: 
\[
\rho (l(T^{n}x-x))\leq k^{n-1}\rho (l(Tx-x))+k^{n-2}r+....+r 
\]
As $\alpha >1$, we have $\rho (l(Tx-x))\leq r$. Then 
\[
\rho (l(T^{n}x-x))\leq \frac{1-k^{n}}{1-k}r 
\]
Therefore, $\rho (c(T^{n+m}x-T^{m}x))\leq k^{m}\frac{1-k^{n}}{1-k}%
r\rightarrow 0$ as $n,m\rightarrow +\infty $\newline
$X_{\rho }$ is $\rho $-complete and $B$ is $\rho $-closed hence, $\exists
z\in B$ such that $\rho (c(T^{n}x-z))\rightarrow 0$ as $n\rightarrow +\infty
.$\newline
We prove that $z$ is a fixed point of $T$. Indeed, 
\begin{eqnarray*}
\rho (\frac{c}{2}(Tz-z)) &=&\rho (\frac{c}{2}(Tz-T^{n+1}x)+\frac{c}{2}%
(T^{n+1}x-z)) \\
&\leq &k\rho (l(z-T^{n}x))+\rho (c(T^{n+1}x-z)) \\
&\leq &\rho (c(z-T^{n}x))+\rho (c(T^{n+1}x-z)) \\
&&
\end{eqnarray*}
since $\rho (c(z-T^{n}x))+\rho (c(T^{n+1}x-z))\rightarrow 0$ as $%
n\rightarrow +\infty $ then $T(\frac{c}{2}(Tz-z))=0$ and $Tz=z$. \newline
Remark 1.1 \newline
It results from this proof that: $\exists x\in B$ such that $\rho
(c(T^{n+m}x-T^{m}x))\leq k^{m}\frac{1-k^{n}}{1-k}r$. If $\rho $ has the
Fatou property, then, the fixed point $z$ is such that: 
\begin{eqnarray*}
\rho (c(z-T^{m}x)) &\leq &\liminf \rho (c(T^{n+m}x-T^{m}x)) \\
&\leq &\frac{k^{m}}{1-k}r
\end{eqnarray*}
This estimate allows an approximation to this fixed point.\newline
Remark 1.2 \newline
If $\rho $ is a $s$-convex modular, we have the same theorem 1.1. But,
because of the $s$-convex combination $(\frac{l}{c})^{s}+\frac{1}{\alpha ^{s}%
}=1$ , some technical modifications are necessary in the theorem 1.1's proof.%
\newline
The comparison between the theorem 1.1 and the theorem I.1 in \cite{ai}
gives the following result.

\begin{coro}
1.1.Let $X_{\rho }$ be a $\rho $-complete modular space, where $\rho $ is $s$%
-convex. $B\subseteq X_{\rho }$ is a \newline
$\rho $-closed subset of $X_{\rho }$. $T:B\rightarrow B$ is a mapping such
that: \newline
$\exists c,k,l\in {R\hspace*{-0.9ex}\rule{0.15ex}{1.5ex}\hspace*{0.9ex}}^{+}$
with $c>Max(l,kl)$ and $\rho (c(Tx-Ty))\leq k^{s}\rho (l(x-y)),\ \ \forall
x,y\in B\ \ \ (\ast \ast )$ \newline
Then $T$ has a fixed point.
\end{coro}

This result brings substantial ameliorations to theorem I.1 \cite{ai}: The
insertion of the constant $l$, the deletion of the hypothesis, the $\Delta
_{2}$-condition and the Fatou property. \newline
Proof of the corollary 1.1: \newline
Let $l_{0}$ be one constant such that $c>l_{0}>Max(l,kl)$; We have: 
\begin{eqnarray*}
\rho (c(Tx-Ty)) &\leq &k^{s}\rho (l(x-y)) \\
&=&k^{s}\rho (\frac{l}{l_{0}}l_{0}(x-y)) \\
&\leq &(\frac{lk}{l_{0}})^{s}\rho (l_{0}(x-y))
\end{eqnarray*}
Then, $c>l_{0}$ and $(\frac{lk}{l_{0}})^{s}<1$ . By the theorem 1.1, $T$ has
a fixed point. \newline
Remark.1.3

It results from the above proof that, if $\rho $ is $s$-convex, the two
formulations of the strong contraction of

$T\ \ (\ (\ast )$ in theorem 1.1. and $(\ast \ast )$ in corollary 1.1.) are
equivalent. \newline
Remark 1.4 \ \ \ \ \ \ \ \ \ \ \ \ \ \ \ \ \ \ \ \ \ \ \ \ \ \ \ \ \ \ \ \ \
\ \ \ \ \ \ \ \ \ \ \ \ \ \ \ \ \ \ \ \ \ \ \ \ \ \ \ \ \ \ \ \ \ \ \ \ \ \
\ \ \ \ \ \ \ \ \ \ \ \ \ \ \ \ \ \ \ \ \ \ \ \ \ \ \ \ \ \ \ \ \ \ \ \ \ \
\ \ \ \ \ \ \ \ \ \ \ \ \ \ \ \ \ \ \ \ \ \ \ \ \ \ \ \ \ \ \ \ \ \ \ \ \ \
\ \ \ \ \ \ \ If $X_{\rho }$ is equipped with the following convergence: $%
x_{n}\stackrel{\rho }{\rightarrow }x\Longleftrightarrow \rho
(x_{n}-x)\rightarrow 0$ as $n\rightarrow +\infty $.\newline
Then the theorem 1.1 takes the following form:

\begin{theo}
1.2.Let $X_{\rho }$ be a $\rho $ complete modular space, and $B\subseteq
X_{\rho }$ a $\rho $-closed subset of $X_{\rho }$. Let $T:B\rightarrow B$ be
a mapping such that: \newline
$\exists c,k,l\in I{R}^{+}{\hspace*{-0.9ex}}$ with $c>l,\ \ k\in ]0,1[$ and $%
\rho (c(Tx-Ty))\leq k\rho (l(x-y));\forall x,y\in B$.\newline
Then $T$ has a fixed point if one of the following assumptions is satisfied
: \newline
i) $1\leq c$ \newline
ii)$0<c<1$ and $\rho $ satisfies the $\Delta _{2}$-condition.
\end{theo}

Proof. \newline
It results from the theorem 1.1.'s proof that $\rho
(c(T^{n+m}x-T^{n}x))\rightarrow 0$ as $n,m\rightarrow +\infty .\ \ X_{\rho }$
is $\rho $-complete, hence, $\exists z\in X_{\rho }$ such that $\rho
(cT^{n}x-z)\rightarrow 0$ as $n\rightarrow +\infty $ . Then $\frac{z}{c}\in
B $. Indeed, $\rho (cT^{n}x-z)=\rho (c(T^{n}x-\frac{z}{c})).$\newline
If $1\leq c$, then $\rho (T^{n}x-\frac{z}{c})\leq \rho (c(T^{n}x-\frac{z}{c}%
))\rightarrow 0$ as $n\rightarrow +\infty $, and $\frac{z}{c}\in B.$ If $%
0<c<1$ and $\rho $ satisfies the $\Delta _{2}$-condition, then $\rho
(c(T^{n}x-\frac{z}{c}))\rightarrow 0$ as $n\rightarrow +\infty \Rightarrow
\rho (T^{n}x-\frac{z}{c})\rightarrow 0$ as $n\rightarrow +\infty $, and $%
\frac{z}{c}\in B$. \newline
We prove that $\frac{z}{c}$ is a fixed point of $T$. Indeed, 
\begin{eqnarray*}
\rho (\frac{c}{2}(T\frac{z}{c}-\frac{z}{c})) &=&\rho (\frac{c}{2}(T\frac{z}{c%
}-T^{n+1}x)+\frac{c}{2}(T^{n+1}x-\frac{z}{c})) \\
&\leq &\rho (c(T\frac{z}{c}-T^{n+1}x))+\rho (cT^{n+1}x-z) \\
&\leq &k\rho (\frac{l}{c}(z-cT^{n}x))+\rho (c(T^{n+1}x-z)) \\
&\leq &k\rho (z-cT^{n}x)+\rho (cT^{n+1}x-z)
\end{eqnarray*}
Since $k\rho (z-cT^{n}x)+\rho (cT^{n+1}x-z)\rightarrow 0$ as $n\rightarrow
+\infty $, then $\rho (\frac{c}{2}(T\frac{z}{c}-\frac{z}{c}))=0$ and $T\frac{%
z}{c}=\frac{z}{c}$.\newline
Remark 4\newline
If $B$ is a subspace of the $X_{\rho }$ in theorem 1.2, then the constraints
on $c$ ($1\leq c$ or $0<c<1$ and $\rho $ satisfies the $\Delta _{2}$%
-condition ) are useless.\newline

\section{II-Strict $\protect\rho $- contraction mappings}

Let us note that if $c=l$ or $c=l=1$ , the adopted method in the theorem
1.1.'s proof is not valid.\newline
The following result can be considered as the second approach of the
Banach's fixed point theorem in the normed spaces. \newline

\begin{theo}
2.1.Let $X_{\rho }$ be a $\rho $-complete modular space where $\rho $
satisfies the $\Delta _{2}$-condition. $B\subseteq X_{\rho }$ a $\rho $%
-closed subset of $X_{\rho }$ . Let $T:B\rightarrow B$ be a strict $\rho $%
-contraction mapping, i.e., \newline
$\exists c,k\in {R\hspace*{-0.9ex}\rule{0.15ex}{1.5ex}\hspace*{0.9ex}}^{+}$
with $k\in ]0,1[$ and $\rho (c(Tx-Ty))\leq k\rho (c(x-y))\ ,\ \forall x,y\in
B$ \newline
We suppose that $\rho (c(x-y))<+\infty ,\ \ \forall x,y\in B.$ Then $T$ has
a unique fixed point.
\end{theo}

Remarks. \newline
This result, As the theorem 1.1, has permitted to delete the boundedness
conditions concerning the domain of $T$ in \cite{krk}-\cite{m1}-\cite{Lm}.%
\newline
But the $\Delta _{2}$-condition, in the following proof, is essential, that
is , the iterative techniques are to happen locally. \newline
Proof of the theorem 2.1 \newline
$1^{st}$ step \newline
If $\rho $ satisfies the $\Delta _{2}$-condition, then, $\exists $ $\delta
,\ \ L,\ \ M\in {\ R\hspace*{-0.9ex}\rule{0.15ex}{1.5ex}\hspace*{0.9ex}}^{+}$
such that: \newline
$\rho (x)\leq \delta \Rightarrow \rho (2x)\leq L\rho (x)+M\ \ \ \ (\Delta
_{2})$ \newline
Otherwise, for $\delta =\frac{1}{n}$ and $L=M=1$, we have $\rho
(x_{n})\rightarrow 0$ as $n\rightarrow +\infty $ and \newline
$\rho (2x_{n})>\rho (x_{n})+1\geq 1$. Absurd.\newline
$2^{nd}$ step\newline
$\exists x_{0}\in B$ such that $r=\rho (2c(Tx_{0}-x_{0}))$ is arbitrary
small, because:\newline
$\rho (c(Tx-x))<\infty ,\rho (c(T^{n+1}x-T^{n}x))\leq k^{n}\rho
(c(Tx-x))\rightarrow 0$ as $n\rightarrow +\infty $ , and, by the $\Delta
_{2} $-condition, $\rho (2c(T^{n+1}-T^{n}x))\rightarrow 0$ as $n\rightarrow
+\infty $. \newline
We suppose that the constant $k$ is such that: 
\[
0<k\leq \frac{\delta }{M+r+L\delta }\ \ \ (1) 
\]
Let us note that $(1)\Rightarrow Lk<1$ . We prove that $\{T^{n}x_{0}\}_{n\!%
\text{ }\in \!N}$ is $\rho $-Cauchy in $X_{\rho }$. Indeed, by induction we
have: 
\[
\rho (c(T^{n+m}x_{0}-T^{m}x_{0}))\leq k^{m}\rho (c(T^{n}x_{0}-x_{0})) 
\]
We show that $\rho (c(T^{n}x_{0}-x_{0}))\leq \frac{1-(Lk)^{n}}{1-Lk}(M+r)\ \
(2)$ \newline
Indeed, for $n=1,\rho (c(Tx_{0}-x_{0}))\leq r\leq M+r$. We suppose that $(2)$
is satified. Then 
\[
\rho (c(T^{n+1}x_{0}-x_{0}))\leq \rho (2c(T^{n+1}x_{0}-Tx_{0}))+r 
\]
Or 
\begin{eqnarray*}
\rho (c(T^{n+1}x_{0}-Tx_{0})) &\leq &k\rho (c(T^{n}x_{0}-x_{0})) \\
&\leq &k\frac{1-(Lk)^{n}}{1-Lk}(M+r) \\
&\leq &\frac{k}{1-Lk}(M+r)
\end{eqnarray*}
By $(1)\ \ ,\ \ \rho (c(T^{n+1}x_{0}-Tx_{0}))\leq \delta $ . Therefore 
\begin{eqnarray*}
\rho (c(T^{n+1}x_{0}-x_{0})) &\leq &Lk\frac{1-(Lk)^{n}}{1-Lk}(M+r)+M+r \\
&\leq &\frac{1-(Lk)^{n}}{1-Lk}(M+r)
\end{eqnarray*}
Finally, we have $\rho (c(T^{n+m}x_{0}-T^{m}x_{0}))\leq k^{m}\frac{1-(Lk)^{n}%
}{1-Lk}(M+r)\rightarrow 0$ as $n,m\rightarrow +\infty $ . \newline
As $X_{\rho }$ is $\rho $-complete and $B$ is $\rho $-closed, $\exists z\in
B $ such that $\rho (c(T^{n}x_{0}-z))\rightarrow 0$ as $n\rightarrow +\infty 
$. Then $z$ is a fixed point of $T$ . Indeed, 
\begin{eqnarray*}
\rho (\frac{c(Tz-z)}{2}) &=&\rho (\frac{c}{2}%
(Tz-T^{n+1}x_{0}+T^{n+1}x_{0}-z)) \\
&\leq &k\rho (c(z-T^{n}x_{0}))+\rho (c(T^{n+1}x_{0}-z))
\end{eqnarray*}
So $k\rho (c(z-T^{n}x_{0}))+\rho (c(T^{n+1}x_{0}-z))\rightarrow 0$ as $%
n\rightarrow +\infty $. Hence, $\rho (\frac{c(Tz-z)}{2})=0$ and $Tz=z$.
Since $\rho (c(x-y))<\infty ,\ \ \forall x,y\in B$, then $z$ is a unique 
\newline
$3^{rd}$ step \newline
As $k^{n}\rightarrow 0$ as $n\rightarrow +\infty $, then, $\exists p_{0}\in
I\!\!N$ such that 
\[
k^{p_{0}}\leq \frac{\delta }{M+r+L\delta } 
\]
We take $S=T^{p_{0}}$ and $k_{0}=k^{p_{0}}$ we have: 
\[
\rho (c(Sx-Sy))\leq k_{0}\rho (c(x-y)),\ \ \forall x,y\in B 
\]
By the same approachs as in the $2^{nd}$ step, we verify that $S$ has a
unique fixed point $z$. Therefore $z$ is also a unique fixed point of $T$. 
\newline
Remark\newline
Let $(\Omega ,\Sigma ,\mu )$ be a measure space. $(L^{\varphi },\rho )$ is
the Musielak-Orlicz space where $\mu $ is $\sigma $-finite and atomless, and 
$\varphi $ is locally integrable. If $\rho $ satisfies the $\Delta _{2}$%
-condition, then, by ([6], theorem 8.14 ), $\exists L,M\in {\
R\hspace*{-0.9ex}\rule{0.15ex}{1.5ex}\hspace*{0.9ex}}^{+}$ such that $\rho
(2x)\leq L\rho (x)+M,\ \ \forall x\in L\varphi .$ \newline
In this case, the constant $\delta $, in the above proof, is arbitrary;
hence, this proof is valid with the constraint: $\exists p_{0}\in I\!\!N$
such that $k^{p_{0}}L<1$

\section{III-$\protect\rho $-nonexpansive mappings}

In this paragraph, we consider the modular space $X_\rho $ equipped with the
convergence:\newline
$x_n \stackrel{\rho}{\rightarrow} x \Longleftrightarrow \rho (x_n -x )\to 0 $
as $n\to +\infty $

\begin{deff}
. The modular $\rho $ satisfies the regular growth condition if $W_{\rho
}(t)<1$ for all $t\in \lbrack 0,1[$, where $W_{\rho }(t)=\sup \{\frac{\rho
(tx)}{\rho (x)},\ \ x\in X_{\rho },\ \ 0<\rho (x)<\infty \}$. \newline
All $s$-convex function modulars satisfy the regular growth condition. For
other examples see \cite{krk} \newline
. The set $B$ is said to be star-shaped if there exists $z\in B$ such that $%
\alpha z+\beta x\in B,\ \ \forall x\in B$, whenever $\alpha ,\beta \in I{%
R\hspace*{-0.9ex}}^{+}$ with $\alpha +\beta =1$. Such a point $z$ is called
a center of $B$ .\newline
A subset $B$ of $X_{\rho }$ is said to be $\rho $-bounded in the sense of
topological vector spaces ($\tau _{\rho }$-bounded) if: For every sequence $%
\{x_{n}\}\subset B$ and any sequence of numbers $\epsilon _{n}\rightarrow 0$%
, there holds $\rho (\epsilon _{n}x_{n})\rightarrow 0$ as $n\rightarrow
+\infty $.
\end{deff}

The following result can be considered as of the Schauder's type.

\begin{theo}
3.1.Let $X_{\rho }$ be a $\rho $-complete modular space where $\rho $
satisfies the regular growth condition. $B\subseteq X_{\rho }$ a $\rho $%
-closed and star-shaped subset of $X_{\rho }$. \newline
$T:B\rightarrow B$ be a $\rho $-nonexpansive mapping, i.e., $\rho
(Tx-Ty)\leq \rho (x-y),\ \ \forall x,y\in B$. \newline
If ${\overline{T(B)}}^{\rho }$ is $\rho $-compact, then $T$ has a fixed
point.
\end{theo}

Remark 3.1\newline
This result is presented under supplementary conditions in \cite{krk}-\cite
{m1}, where $X_{\rho }=L_{\rho }$; in particular, $B$ is $\rho $-bounded ( $%
\delta _{\rho }(B)={\displaystyle\sup_{x,y\in B}}\rho (x-y)<\infty $) and $%
\rho $-compact. We note that if $B$ is $\rho $-compact, then $B$ is $\rho $%
-closed and $T(B)$ is $\rho $-compact. \newline
Proof of the theorem 3.1\newline
$1^{st}$ step

\begin{lem}
3.1.$X_{\rho }$ and $B$ are as in theorem 3.1. $T:B\rightarrow B$ is $\rho $%
-nonexpansive. Then\newline
a) The equation $x=\alpha z+\beta Tx\ \ (z$ a center of $B,\ \ x\in B,\ \
(\alpha ,\beta )\in I{R\hspace*{-0.9ex}\hspace*{0.9ex}}^{+}\times I{%
R\hspace*{-0.9ex}\hspace*{0.9ex}}^{+}$ with $\alpha +\beta =1)$ has one
solution.\newline
b) If moreover, $TB$ is $\tau _{\rho }$-bounded, then $T$ has an
approximating fixed point.
\end{lem}

Proof\newline
a) For $\beta \in ]0,1[$, let $\lambda \in ]1,\frac{1}{\beta }[$. We
consider $Sx=\alpha z+\beta Tx$. Then $S:B\rightarrow B$ and 
\begin{eqnarray*}
\rho (\lambda (Sx-Sy)) &=&\rho (\lambda \beta (Tx-Ty)) \\
&\leq &W_{\rho }(\lambda \beta )\rho (Tx-Ty) \\
&\leq &W_{\rho }(\lambda \beta )\rho (x-y)
\end{eqnarray*}
By the theorem 1.2, $S$ has a fixed point \newline
b) Let $k_{n}\in ]0,1[$ with $k_{n}\nearrow 1$. By a), we have $%
x_{n}=(1-k_{n})z+k_{n}Tx_{n}$. Hence \newline
$\rho (Tx_{n}-x_{n})\leq \rho (2(1-k_{n})Tx_{n})+\rho (2(1-k_{n})z)$ \newline
By the definition of $X_{\rho }$, $\rho (2(1-k_{n})z)\rightarrow 0$ as $%
n\rightarrow +\infty $. As $TB$ is $\tau _{\rho }$-bounded, then, \newline
$\rho (2(1-k_{n})Tx_{n})\rightarrow 0$ as $n\rightarrow +\infty $. \newline
Therefore $\rho (Tx_{n}-x_{n})\rightarrow 0$ as $n\rightarrow +\infty $ i.e. 
$T$ has an approximating fixed point . \newline
$2^{nd}$ step \newline
As ${\overline{T(B)}}^{\rho }$ is $\rho $-compact, then ${\overline{T(B)}}%
^{\rho }$ is $\tau _{\rho }$-bounded. [jm]. Hence, there exists $x_{n}\in B$
such that $\rho (Tx_{n}-x_{n}))\rightarrow 0$ as $n\rightarrow +\infty .$
Also, there exists $\{Tx_{n^{\prime }}\}$ a subsequence of $\{Tx_{n}\}$ that
is $\rho $-convergent to $y\in B$. We prove that $y$ is a fixed point of $T$%
. Indeed, we have: 
\begin{eqnarray*}
\rho (\frac{Ty-y}{3}) &=&\rho (\frac{(Ty-T^{2}x_{n^{\prime
}})+(T^{2}x_{n^{\prime }}-Tx_{n^{\prime }})+(Tx_{n^{\prime }}-y)}{3}) \\
&\leq &2\rho (Tx_{n^{\prime }}-y)+\rho (Tx_{n^{\prime }}-x_{n^{\prime }})
\end{eqnarray*}
Since $2\rho (Tx_{n^{\prime }}-y)+\rho (Tx_{n^{\prime }}-x_{n^{\prime
}})\rightarrow 0$ as $n^{\prime }\rightarrow +\infty $ . Therefore $\rho (%
\frac{Ty-y}{3})=0$ and $Ty=y$ .\newline
Remark 3.3\newline
Recall that $B$ is $\rho $-bounded is, in general, not equivalent to $B$ is $%
\tau _{\rho }$-bounded, see \cite{Lm}.\newline
Finally, we present one result, for the strict $\rho $-contration mappings,
using the lemma 3.2. \newline

\begin{prop}
3.1.Let $X_{\rho }$ be a $\rho $-complete modular space where $\rho $ is
convex subset of $X_{\rho }$ with $0\in B$ . \newline
Let $T:B\rightarrow B$ be a mapping such that $\exists $ $k\in ]0,1[$ with $%
\rho (Tx-Ty)\leq k\rho (x-y);\ \ \forall x,y\in B$\newline
Let ${\cal {A}}$ $=\{x\in B:x=\lambda Tx,\ \ \lambda \in ]0,1[\}$. If ${%
\displaystyle\sup_{x\in {\cal {A}}}}\rho (x)<\infty $, then $T$ has a fixed
point.
\end{prop}

Proof. \newline
By the lemma 3.2, ${\cal {A}\neq \phi }$. Let $\lambda _{n}\in ]0,1[$, with $%
\lambda _{n}\nearrow 1$, and $x_{n}=\lambda _{n}Tx_{n}$. We show that $%
\{x_{n}\}$ is $\rho $-Cauchy. Indeed, for $m>n$, we have 
\begin{eqnarray*}
\rho (x_{m}-x_{n}) &=&\rho (\lambda _{m}Tx_{m}-\lambda _{n}Tx_{n}) \\
&=&\rho (\lambda _{n}(Tx_{m}-Tx_{n})+(\lambda _{m}-\lambda _{n})Tx_{m}) \\
&=&\rho (\lambda _{n}(Tx_{m}-Tx_{n})+\frac{\lambda _{m}-\lambda _{n}}{%
\lambda _{m}}x_{m})
\end{eqnarray*}
As $\lambda _{n}+\frac{\lambda _{m}-\lambda _{n}}{\lambda _{m}}\leq 1$, then 
\[
\rho (x_{m}-x_{n})\leq \lambda _{n}k\rho (x_{m}-x_{n})+\frac{\lambda
_{m}-\lambda _{n}}{\lambda _{m}}\rho (x_{m}) 
\]
So $\rho (x_{m}-x_{n})\leq \frac{\lambda _{m}-\lambda _{n}}{\lambda _{m}(1-k)%
}{\sup_{m}}\rho (x_{m})\rightarrow 0$ as $m,n\rightarrow +\infty $ . \newline
$X_{\rho }$ is a $\rho $-complete space and $B$ is a $\rho $-closed, hence, $%
\exists x\in B$ such that $\rho (x_{n}-x)\rightarrow 0$ as $n\rightarrow
+\infty $. We show that $x$ is a fixed point of $T$. Indeed, we have 
\begin{eqnarray*}
\rho (\frac{Tx-x}{3}) &=&\rho (\frac{Tx-Tx_{n}+Tx_{n}-x_{n}+x_{n}-x}{3}) \\
&\leq &(k+1)\rho (x-x_{n})+\rho (Tx_{n}-x_{n})
\end{eqnarray*}
Or $Tx_{n}-x_{n}=x_{n}(\frac{1}{\lambda _{n}}-1)$ and for $n$ very large, we
have $0<\frac{1}{\lambda _{n}}-1<1$; hence \newline
$\rho (Tx_{n}-x_{n})\leq (\frac{1}{\lambda _{n}}-1)\ {\sup_{n}}\rho
(x_{n})\rightarrow 0$ as $n\rightarrow +\infty $ . Therefore $\rho (\frac{%
Tx-x}{3})=0$ and $Tx=x$.

\ \ \ \ \ \ \ \ \ \ \ \ \ \ 

\ \ \ \ \ \ \ \ \ \ \ \ E-mail: hanebaly@hotmail.com\ . \ \ \ 

\ \ \ \ \ \ \ \ \ \ \ 

\ \ \ \ \ \ \ \ \ \ \ \ \ \ \ \ \ \ \ \ \ \ \ \ \ \ \ \ \ \ \ \ \ \ \ \ \ \
\ \ \ \ \ \ \ \ \ \ \ \ \ \ \ \ \ \ \ \ \ \ \ \ \ \ \ \ \ \ \ 

\end{document}